\documentclass[11pt, reqno]{amsart}

\usepackage{amssymb}
\usepackage{hyperref} 
\usepackage{mathrsfs,bbold}
\usepackage{amsmath}
\usepackage{amsthm}
\usepackage{enumerate}
\usepackage{dsfont}
\usepackage{color}
\usepackage[a4paper]{geometry}
\usepackage[utf8]{inputenc}
\usepackage{stmaryrd}
\usepackage{titlesec}
\usepackage{mathabx}
\usepackage{appendix}
\usepackage{todonotes, fancyhdr}

\usepackage[all]{xy}
\let\oldtocsection=\tocsection
\let\oldtocsubsection=\tocsubsection
\let\oldtocsubsubsection=\tocsubsubsection

\renewcommand{\tocsection}[2]{\hspace{0em}\oldtocsection{#1}{#2}}
\renewcommand{\tocsubsection}[2]{\hspace{1em}\oldtocsubsection{#1}{#2}}
\renewcommand{\tocsubsubsection}[2]{\hspace{2em}\oldtocsubsubsection{#1}{#2}}

\numberwithin{theorem}{section}
\numberwithin{equation}{section}

\newcommand{\ssk}{\smallskip}

\renewcommand{\epsilon}{\varepsilon}

\newcommand\bbR{\mathbb{R}}

\newcommand{\mcB}{\mathcal{B}}

\newcommand{\mcF}{\mathcal{F}}
\newcommand{\mcG}{\mathcal{G}}

\newcommand\mcT{\mathcal T}

\newcommand{\bfX}{\textsf{\textbf{X}}}

\newenvironment{Dem}[1][\unskip]{%
    \begin{list}{\hspace{0.5cm}{\sf \textbf{Proof #1 --}}}{%
        \setlength{\topsep}{0pt}%
        \setlength{\leftmargin}{0pt}%
        \setlength{\rightmargin}{0pt}%
        \setlength{\listparindent}{0pt}%
        \setlength{\itemindent}{0pt}%
        \setlength{\parsep}{0pt}%
        \addtolength{\leftmargin}{20pt}%
        \addtolength{\rightmargin}{0pt}%
    } \item }{\hfill $\rhd$\end{list}\smallskip}

\titleformat{\section}[block]
{\filcenter\normalfont\sffamily\bfseries\Large}
{{\hspace{-0.7cm}}\thesection \hspace{0.2em} --\vspace{0.3cm}}{0.5em}{}

\titleformat{\subsection}[block]
{\filcenter\normalfont\sffamily\bfseries\large}  						  
{\hspace{-0.7cm}\thesubsection \hspace{0.5em}--\vspace{0.3cm}}{.5em}{}  
\titlespacing{\subsection}{-0pc}{1.5ex plus .1ex minus .2ex}{0pc}

\titleformat{\subsubsection}[block]
{\normalfont\sffamily\bfseries}{\hspace{-1cm}\thesubsubsection}{.5em}{}
\titlespacing{\subsubsection}{15pc}{1ex plus .1ex minus .2ex}{1pc}

\def\XXint#1#2#3{{\setbox0=\hbox{$#1{#2#3}{\int}$}
     \vcenter{\hbox{$#2#3$}}\kern-.5\wd0}}

\numberwithin{subsection}{section}
\numberwithin{subsubsection}{subsection}

\newtheoremstyle{mystyle}
{3pt}               
{3pt}               
{\it }                      
{}                      
{\sffamily\bfseries}             
{}                      
{0.5em}                 
{\llap{#2. }#1{$\;$ --}}

\theoremstyle{mystyle}

\newtheorem{thm}{Theorem}
\newtheorem*{thm*}{Theorem}

\newtheorem{lem}[thm]{\hspace{-0.2cm}  {Lemma} }
\newtheorem*{defn*} {Definition}
\newtheorem*{prop*} {Proposition}
\newtheorem*{lem*} {Lemma}
\newtheorem*{cor*} {Corollary}

\newtheoremstyle{mystyle2}
{3pt}               
{3pt}               
{\it }                      
{}                      
{\sffamily\bfseries}             
{}                      
{0.5em}                 
{\llap{#2 }#1{\hspace{0.2cm}--}}
\theoremstyle{mystyle2}

\newtheorem*{definition*}{Definition}
\newtheorem*{theorem*}{Theorem}
\numberwithin{equation}{section} 
\newtheorem*{Remark*}{Remark}

\newcommand\restr[2]{{
  \left.\kern-\nulldelimiterspace 
  #1 
  \vphantom{\big|} 
  \right|_{#2} 
  }}

\newcommand{\bR}{\mathbf{R}}

\newcommand{\RR}{\mathbb{R}}

\begin{document}

\begin{center}
{\Huge\sffamily{On the definition of a solution to a rough differential equation   \vspace{0.5cm}}}
\end{center}
\vskip 5ex minus 1ex

\begin{center}
{\sf I. BAILLEUL}\footnote{I.Bailleul thanks the Centre Henri Lebesgue ANR-11-LABX-0020-01 for its stimulating mathematical research programs, and the U.B.O. for their hospitality, part of this work was written there. Partial support from the ANR-16-CE40-0020-01.}\footnote{Univ. Rennes, CNRS, IRMAR - UMR 6625, F-35000 Rennes, France.}
\end{center}

\vspace{1cm}

\begin{center}
\begin{minipage}{0.8\textwidth}
\renewcommand\baselinestretch{0.7} \scriptsize \textbf{\textsf{\noindent Abstract.}} 
There are several ways of defining what it means for a path to solve a rough differential equation. The most commonly used notion is due to Davie; it involves a Taylor expansion property that only makes sense a priori in a given coordinate system. Bailleul's definition \cite{BailleulRMI} is coordinate independent. Cass and Weidner \cite{CassWeidner} recently proved that the two definitions are actually equivalent, using deep algebraic insights on rough paths. We provide in this note an algebraic-free elementary short proof of this fact.
\end{minipage}
\end{center}

\vspace{0.5cm}

\begin{center}
\begin{minipage}{0.8\textwidth}
\renewcommand\baselinestretch{0.7} \scriptsize \textbf{\textsf{\noindent R\'esum\'e.}} 
Diff\'erentes notions de solution d'une \'equations diff\'erentielle rugueuses sont disponibles dans la litt\'erature. La notion la plus utilis\'ee est due \`a Davie ; elle met en jeu un d\'eveloppement de Taylor, dont le sens est li\'e \`a un choix de coordonn\'ees. La d\'efinission donn\'ee par Bailleul dans \cite{BailleulRMI} ne d\'epend pas d'un choix de coordonn\'ees. Cass et Weidner \cite{CassWeidner} ont r\'ecemment d\'emontr\'e que les deux d\'efinitions sont en fait \'equivalentes, en s'appuyant sur des r\'esultats alg\'ebriques \'elabor\'es. On donne dans cette note une d\'emonstration courte et \'el\'ementaire de ce fait.

\end{minipage}
\end{center}

\vspace{1cm}

\section{Main result}
\label{SectionMain}

\textit{a) Setting -- }Rough paths theory was introduced by T. Lyons in \cite{Lyons98} as a theory of controlled differential equations 
\begin{equation}
\label{EqRDE}
\dot x_t = \textrm{F}(x_t)\dot h_t
\end{equation}
with $\textrm{F}\in L(\RR^\ell,\RR^d)$, and a non-differentiable control $h$ with values in $\RR^\ell$. The path $x$ takes here its values in $\RR^d$. Variants of the Cauchy-Lipschitz theorem gives the well-posed character of this equation for an absolutely continuous control, and L.C. Young extended this analysis to controls that are $\alpha$-H\"older, for $\alpha>\frac{1}{2}$; the solution path is then a continuous function of the control $h\in C^{\alpha}([0,1],\RR^\ell)$, under proper regularity and boundedness conditions on F. Nothing better than that can be done in a deterministic setting in the H\"older class, on account of the fact that since the solution path $x_t$ is expected to be no better than $\alpha$-H\"older, so is $\textrm{F}(x_t)$, so the equation involves making sense of the product $\textrm{F}(x_t)\dot h_t$ of an $\alpha$-H\"older function with an $(\alpha-1)$-H\"older distribution, either directly or in its integral form. It is known that no such product can be defined as a continuous function of its two arguments when the sum of the H\"older regularity exponents add up to a non-positive real number -- see e.g. \cite{BCD}. Lyons' deep insight was to realize that what really controls the dynamics in \eqref{EqRDE} is not the path $h$, but rather the data of the path together with a number of its iterated integrals 
$$
\int h^i_sdh^j_s, \quad \iint_{r\leq s} h^i_rdh^j_rdh^k_s, \dots
$$
for all possible indices $1\leq i,j,k,\dots\leq \ell$. The rougher the control $h$ is, the more iterated integrals you need. This can be understood from the fact that these iterated integrals are the coefficients that appear in a Taylor expansion of the solution to equation \eqref{EqRDE} when the control is smooth, and this leads to the notion of solution adopted in \cite{Davie, BailleulRMI} by Davie and Bailleul. There is no miracle though, and one cannot still make sense of these iterated integrals when the control is not sufficiently regular. Lyons' considerable feat was to extract from these 'non-existing' quantities the analytic and algebraic properties that they should satisfy and to work directly with objects enjoying these properties as new controls. These are the weak geometric H\"older $p$-rough paths that we encounter below. We refer the reader to the gentle introductions \cite{LyonsStFlour, FrizHairer, LNRoughPaths} for different points of views on rough differential equations; one can learn everything from scratch in the first 19 pages of \cite{BailleulRMI}. 

\medskip

\textit{b) Weak geometric H{\"o}lder $p$-rough paths and rough differential equations --} We refer the reader to \cite{Lyons98, BaudoinLN, BailleulRMI} for basics on rough paths and recall all we need in this section. Pick an integer $p\geq 1$. Denote by $(e_i)_{i=1..\ell}$ the canonical basis of $\bbR^\ell$, and by $T^{[p]}_\ell=\bigoplus_{i=1}^\ell(\bbR^\ell)^{\otimes i}$ the truncated tensor algebra over $\bbR^\ell$; it is equipped with the Lie bracket operation $[a,b]:=ab-ba$, and endowed with a norm 
$$
\|a\| := \sum_{i=1}^{[p]} |a^i|^{1/i},
$$
for $a=:\oplus_{i=1}^{[p]} a^i$; denote by $(a^I)_{I\in\llbracket 1,\ell\rrbracket^k, 0\leq k\leq [p]}$ the coordinates of a generic element $a$ of $T^{[p]}_\ell$ in the canonical basis of $T^{[p]}_\ell$. We sometimes write $a=(a^i)_{i=1..[p]}$ instead of $a=\oplus_{i=1}^{[p]} a^i$. Consider $\bbR^\ell$ as a subset of $T^{[p]}_\ell$. The Lie algebra $\frak{g}^{[p]}_\ell$ generated by $\bbR^\ell$ within $T^{[p]}_\ell$ is called the $p$-step free nilpotent Lie algebra, and its exponential $G^{[p]}_\ell$ is called the $p$-step free nilpotent Lie group. Recall that a weak geometric H\"older $p$-rough path $\bfX$ over $\RR^\ell$ is a $\frac{1}{p}$-H\"older path with values in $G^{[p]}_\ell$. Set ${\bfX}_{ts}:={\bfX}_s^{-1}{\bfX}_t$, so we have $\|{\bfX}_{ts}\| \leq O\big(|t-s|^{1/p}\big)$; this is equivalent to having 
$$
\|{\bfX}^i_{ts}\| \leq O\big(|t-s|^{i/p}\big).
$$ 
The archetype of a weak geometric H\"older $p$-rough path is given by the canonical lift of a smooth $\bbR^\ell$-valued path $h$ under the form 
\begin{equation}
\label{EqLiftSmooth}
H_t := \left(1,h_t-h_0,\int_0^t (h_{u_1}-h_0)\,du_1,\dots,\int_{0\leq u_{[p]}\leq \cdots\leq u_1\leq t}dh_{u_{[p]}}\cdots dh_{u_1}\right);
\end{equation}
it can be seen, after Chen \cite{Chen}, that $H$ takes values in $G^{[p]}_\ell$ and that 
$$
H_{ts} := H_s^{-1}H_t = \left(1,h_t-h_s,\int_s^t (h_{u_1}-h_0)\,du_1,\dots,\int_{s\leq u_{[p]}\leq \cdots\leq u_1\leq t}dh_{u_{[p]}}\cdots dh_{u_1}\right).
$$
Denote by $\llbracket1,\ell\rrbracket$ the set of integers between $1$ and $\ell$. Define inductively on the size of the tuple $J := \{i_1,I\}\in\llbracket 1,\ell\rrbracket^{k+1}$, the element
$$
e_{[J]} = \big[e_{i_1},e_{[I]}\big]
$$
of $\frak{g}^{[p]}_\ell$. Given $0\leq s\leq t\leq T$, set 
$$
{\bf \Lambda}_{ts} := \log {\bf X}_{ts} = \underset{0\leq k\leq [p]}{\sum_{I\in\llbracket 1,\ell\rrbracket^k}} \Lambda^I_{ts} e_{[I]} \; \in \frak{g}^{[p]}_\ell.
$$

Let $\textrm{F}=\big(V_1,\dots, V_\ell\big)$ stand for a collection of $\ell$ smooth enough vector fields on $\RR^d$, and $\bfX$ stand for a weak geometric H\"older $p$-rough path  over $\RR^\ell$. There are several definitions of a solution to the rough differential equation 
\begin{equation}
\label{EqRDE}
dz_t = \textrm{F}(z_t)\,d{\bfX}_t.
\end{equation}
Besides Lyons' original definition \cite{Lyons98} and Gubinelli's formulation in terms of controlled paths \cite{Gubinelli}, Davie gave in \cite{Davie} a formulation of a solution in terms of Taylor expansion for the time increment of a solution. We identify freely a vector field to a first order differential operator via the expression
$$
V_if = (Df)(V_i)
$$
for any $C^1$ function $f$ on $\RR^d$. For a tuple $I = (i_1,\dots,i_k)\in\llbracket 1,\ell\rrbracket^k$, we write $\vert I\vert := k$, and define a $k^\textrm{th}$ order differential operator setting
$$
V_I = V_{i_1}\cdots V_{i_k},
$$
that is to say $V_If = V_{i_1}\big(\cdots (V_{i_k}f)\big)$. With a slight abuse of notation, we write $V_I(x)$ for $\big(V_I\textrm{Id}\big)(x)$. 
An $\RR^d$-valued path $(z_t)_{0\leq t\leq T}$ defined on some finite time interval $[0,T]$ is a solution to the rough differential equation \eqref{EqRDE} in the sense of Davie if one has
\begin{equation}
\label{EqTaylorDefn}
z_t = z_s + \sum_{I;\vert I\vert \leq [p]} X^I_{ts}(V_I)(z_s) + O\big(\vert t-s\vert^a\big),
\end{equation}
for any $0\leq s\leq t\leq T$, for some exponent $a>1$. The function $O(\cdot)$ is allowed to depend on F. Two paths satisfying that condition with different exponents $a,a'>1$ coincide if they start from the same point. This definition is a priori not well suited to make sense of a solution to a rough differential taking values in a manifold -- an expression like $m+V(m)$, for a point $m$ of a manifold $M$ and a vector $V(m)\in T_mM$, has for instance no intrinsic meaning. Bailleul gave in \cite{BailleulRMI} a more general notion of solution to a rough differential equation by requiring from a potential solution that it satisfies the estimate
$$
f(z_t) = f(z_s) + \sum_{I;\vert I\vert \leq [p]} X^I_{ts}(V_If)(z_s) + O_f\big(\vert t-s\vert^a\big),
$$
for any real-valued sufficiently regular function $f$ defined on the state space. The function $O_f(\cdot)$ is allowed to depend on $f$. The exponent $a$ that appears here may be different from the exponent that appears in equation \eqref{EqTaylorDefn}; this makes no difference as long as $a>1$, as noted above. This definition makes perfect sense in a manifold setting, and taking $f$ to be a coordinate function in the above equation makes it is clear that a solution in this sense to an equation with values in $\RR^d$ is also a solution in the sense of Davie. 

\medskip

\textit{c) Main result --} Cass and Weidner proved in Theorem 5.3 of their recent work \cite{CassWeidner} that the two notions are actually equivalent using some deep insights from algebra based on the use of Grossman-Larson and Connes-Kreimer Hopf algebras. That point can be proved by elementary means. We assume as usual that the vector fields $V_i$ in \eqref{EqRDE} are $C^\gamma$, for some $\gamma>p$.

\medskip

\begin{thm}
\label{Thm}
An $\bR^d$-valued path $z$ is a solution to the rough differential equation \eqref{EqRDE} in the sense of Davie if and only if it is a solution of that equation in the sense of Bailleul.
\end{thm}

\medskip

\begin{Dem}
We only need to provee that solutions to equation \eqref{EqRDE} in Davie' sense are solutions in Bailleul' sense. Given a globally Lipschitz vector field $V$ on $\RR^d$, we denote by $\exp(V)$ the time $1$ map of the differential equation 
$$
\dot y_u = V(y_u),
$$
that associates to $x$ the value at time $1$ of the solution to the equation started from $x$. Also, define inductively on the size of the tuple $\{i_1,I\}\in\llbracket 1,\ell\rrbracket^{k+1}$, the vector fields
$$
V_{[\{i_1,I\}]} := \big[V_{i_1},V_{[I]}\big],
$$
starting with $V_{[i]} := V_i$. Given $0\leq s\leq t\leq T$, recall we write ${\bf \Lambda}_{ts}$ for $\log {\bf X}_{ts}$; this is an element of $\frak{g}_\ell^{[p]}$. It is proved in Proposition 9 of \cite{BailleulRMI} that one has
\begin{equation}
\label{EqUpperBound}
\left\| \exp\left(\sum_{I;\vert I\vert \leq [p]} \Lambda^I_{ts} V_{[I]}\right) - \Big\{\textrm{Id} + \sum_{I;\vert I\vert \leq [p]} X^I_{ts}(V_I)\Big\}\right\|_{L^\infty} \leq O\Big(\vert t-s\vert^\frac{\gamma}{p}\Big)
\end{equation}
for a $O(\cdot)$ depending only on the vector fields $V_i$. (See also \cite{BoutaibLyonsYang} for a similar statement.) A path $z$ is thus a solution to the rough differential equation \eqref{EqRDE} in the sense of Davie if and only if 
$$
z_t = \exp\left(\sum_{I;\vert I\vert \leq [p]} \Lambda^I_{ts} V_{[I]}\right)(z_s) + O\big(\vert t-s\vert^{a'}\big),
$$
for some exponent $1<a'\leq \frac{\gamma}{p}$. Let $y$ stand now for a solution path to the ordinary differential equation
$$
\dot y_u = \sum_{I;\vert I\vert \leq [p]} \Lambda^I_{ts} V_{[I]}(y_u),
$$ 
started from $z_s$ at time $0$. We use in the computation below the notation $O\big(\vert t-s\vert^{a'}\big)$ for a function whose value may change from line to line. For a $\gamma$-H\"older real-valued function $f$ on $\bbR^d$ one has 
\begin{equation*}
\begin{split}
f(z_t) &= f(y_1) +  O\big(\vert t-s\vert^{a'}\big)   \\
		  &= f(z_s) + \int_0^1 \frac{d}{du}f(y_u)\,du + O\big(\vert t-s\vert^{a'}\big)   \\
		  &= f(z_s) + \sum_{\vert I_1\vert\leq [p]}\Lambda^{I_1}_{ts}\int_0^1V_{I_1}(y_{u_1})\,du_1 + O\big(\vert t-s\vert^{a'}\big)   \\
		  &= f(z_s) + \sum_{\vert I_1\vert\leq [p]}\Lambda^{I_1}_{ts}V_{I_1}(z_s) + \underset{\vert I_1\vert+\vert I_2\vert\leq [p]}{\sum_{\vert I_1\vert\leq [p], \vert I_2\vert\leq [p]}} \int_0^1V_{I_2}V_{I_1}(y_{u_2})\,du_2du_1 + O\big(\vert t-s\vert^{a'}\big),
\end{split}
\end{equation*}
and, after repeating $[p]$ times the same computation giving a function as the integral of its derivative, one eventually has
\begin{equation*}
\begin{split}
f(z_t) &= f(z_s) + \sum_{k=1}^{[p]}\frac{1}{k!}\,\Lambda^{I_1}_{ts}\cdots \Lambda^{I_k}_{ts}\,\big(V_{I_k}\cdots V_{I_1}\big)(z_s) + O\big(\vert t-s\vert^{a'}\big)   \\
		  &= f(z_s) + \sum_{I;\vert I\vert \leq [p]} X^I_{ts}(V_If)(z_s) + O_f\big(\vert t-s\vert^{a'}\big),
\end{split}
\end{equation*}
from the fact that ${\bf \Lambda}_{ts} = {\bf \log}\bfX_{ts}$. We repeat that computation in Section \ref{SectionBranched} in a more general setting.
\end{Dem}

\bigskip

The following remarks emphasize the robust character of the above proof.   \vspace{0.1cm}

\begin{itemize}
   \item One can weaken the regularity assumptions on the vector fields $V_i$ by only requiring that they are $C^{[\gamma]}$ in the usual sense and the $[\gamma]$-derivative of $V_i$ is $(\gamma-[\gamma])$-H\"older continuous, without requiring that the $V_i$ or their derivatives be bounded. Indeed, given a solution $z$ of \eqref{EqRDE} in the sense of Davie, and a time $s$, the classical Cauchy Lipschitz theory ensures that  $\exp\left(\sum_{I;\vert I\vert \leq [p]} {\bf \Lambda}_{ts} V_{[I]}\right)(z_s)$ is well-defined for $t$ close enough to $s$, possibly depending on $z_s$.   \vspace{0.1cm}

   \item For a non-necessarily H\"older weak geometric $p$-rough path $\bfX$ controlled by a general control $w(s,t)$ rather than by $\vert t-s\vert^{1/p}$, one replaces $(t-s)^{1/p}$ by $w(s,t)$ in the definitions of a solution to a rough differential equation driven by $\bfX$, and requires in both definitions that the remainder is of order $w(s,t)^{ap}$, for some exponent $a>1$. The proof of Theorem \ref{Thm} remains the same. The equivalence problem for the different definitions of a solution to a rough differential equation is further explored in the very recent work \cite{BraultLejay} of Brault and Lejay.   \vspace{0.1cm}

   \item Let $E$ stand for a Banach space. For dynamics driven by an $E$-valued rough paths, following \cite{CDLL} and \cite{CassWeidner}, one needs to define $T^{(n)}(E)$ as the completion of the truncated algebraic tensor algebra with respect to a system of cross (semi)norms. One then defines the $n$-step free Lie algebra over $E$ as the closure in $T^{(n)}(E)$ of the algebraic free Lie algebra. The $n$-step free nilpotent Lie group is then indeed the image by the exponential map of the closed Lie algebra. An appropriate formalism where vector fields are indexed by the elements of $E$ is required, as in \cite{BailleulSeminaire, BoutaibLyonsYang} or \cite{CassWeidner}. Identity \eqref{EqUpperBound} is proved in this setting in \cite{BailleulSeminaire} by elementary means similar to their finite dimensional analogues. (We erronously worked in the non-complete versions of the Lie and truncated tensor algebras in \cite{BailleulSeminaire}. The complete setting should be adopted, and nothing is changed to the story told in \cite{BailleulSeminaire} in this extended setting. Thanks to T. Cass for pointing this out.)   \vspace{0.1cm}
   
   \item One of the nice points of Cass and Weidner's work \cite{CassWeidner} is the fact that they can also handle differential equations driven by a more general notion of rouhg path called branched rough paths. In a finite dimensional setting, one can appeal to Hairer and Kelly's result \cite{HairerKelly} stating that a solution to a rough differential equation driven by a branched rough path over $\RR^\ell$ is also the solution to a differential equation driven by a weak geometric rough path over a larger space $\RR^{\ell'}$, with driving vector fields built from the initial $V_i$. The equivalence of the two notions of solution for equations driven by weak geometric rough paths implies the equivalence of the two notions of solution for the equation driven by a branched rough path, since the Taylor expansion of $f(y_t)$ from the branched rough path and from the weak geometric rough path point of views coincide, by Theorem 5.8 in \cite{HairerKelly}. We give a self-contained proof of the equivalence of a Davie-type and Bailleul-type notion of a solution to a rough differential equation driven by a branched rough path in Section \ref{SectionBranched}.
\end{itemize}

\bigskip

\section{Branched rough paths and their associated flows in a nutshell}
\label{SectionBranched}

We prove in this section the following analogue of Theorem \ref{Thm}.

\begin{thm}
\label{Thm2}
An $\bR^d$-valued path $z$ is a solution to a rough differential equation driven by a branched rough path in the sense of Davie if and only if it is a solution of that equation in the sense of Bailleul.
\end{thm}

Theorem \ref{Thm2} gives an alternative proof of Theorem 5.3 of \cite{CassWeidner} in the setting of dynamics driven by branched rough paths. We include in this section a short self-containted introduction to finite dimensional branched rough paths before proving Theorem \ref{Thm2} in Section \ref{SubsectionProof}. We refer to Gubinelli's original paper \cite{GubinelliBranched} and Hairer and Kelly's paper \cite{HairerKelly} for alternative accounts on branched rough paths in a finite dimensional setting, and to \cite{CassWeidner} for an account of the theory in an infinite dimensional setting. 

\medskip

\textit{What are branched rough paths good for?} Together with Equation \eqref{EqLiftSmooth} giving the lift of a smooth $\bbR^\ell$-valued path to the $[p]$-step free nilpotent Lie group, the definition of a weak geometric rough path $\bfX$ makes it clear that the higher order levels of $\bfX$ play the role of the non-existing integrals $\int X^idX^j, \iint X^idX^jdX^k$, etc. The archetypal example of a non-smooth weak geometric $p$-rough path is given by realizations of the Brownian rough path ${\bf B} = (B,\mathbb{B})$, for $2<p<3$, defined by 
$$
\mathbb{B}_{ts}^{jk} := \int_s^t(B_u-B_s)\,{\circ d}B_u.
$$
Note the use of Stratonovich integration. If one defines instead the iterated integral $\mathbb{B}$ with an It\^o integral, the corresponding object $\bf B$ no longer takes values in the $[p]$-step free nilpotent Lie group $G^{[p]}_\ell$. This is related to the fact that It\^o integration procedure does not satisfy the usual chain rule. A different algebraic setting is needed to handle the integration theory of such controls; this is precisely what branched rough paths provide us with. The name ''branched rough paths'' come from the fact that these objects are indexed by trees.

\medskip

\subsection{Branched rough paths}

\textbf{\textsf{1. The starting point --}} Consider an ordinary controlled $\RR^d$-valued differential equation 
$$
\dot x_r = \sum_{i=1}^\ell V_i(x_r)\dot h^i_r
$$
driven by some smooth $\RR^\ell$-valued path $h$. Given any tuple $I\in\llbracket1,\ell\rrbracket^k$, write $ H^I_{ts}$ for the iterated integral $\int_{s\leq s_1\leq \cdots\leq s_k\leq t}dh^{i_1}_{s_1}\dots dh^{i_k}_{s_k}$. Branched rough path appears naturally if one expands in the Taylor formula 
\begin{equation}
\label{EqExpansion}
f(x_t) = f(x_s) + \sum_{|I|\leq k_0}  H^I_{ts} (V_If)(x_s) + O\big(|t-s|^{k_0+1}\big)
\end{equation}
all the iterated derivatives $V_If$ in terms of the derivatives of its different terms. We have for instance 
\begin{equation}
\label{EqDerivative}
V_1\big(V_2f\big)(x) = (D_xf)\big((V_1V_2)(x)\big) + (D_x^2f)\big(V_1(x),V_2(x)\big).
\end{equation}
The best way to represent the resulting sum in \eqref{EqExpansion} is to index it by labelled rooted trees. Let $\mcT$ stand for the set of possibly empty rooted labelled trees, with labels in $\{1,\dots,\ell\}$. Given some trees $\tau_1, \dots, \tau_k$ in $\mcT$ and $a\in\{1,\dots,\ell\}$, denote by $[\tau_1\dots\tau_k]_a$ the element of $\mcT$ obtained by attaching the trees $\tau_1, \dots, \tau_k$ to a new root with label $a$. Any element of $\mcT$ can be constructed in this way starting from the empty tree ${\bf 1}$. We denote by $|\tau|$ the number of vertices of a tree $\tau$, and define its symmetry factor $\sigma(\tau)$ recursively
$$
\sigma\big([\tau_1^{n_1}\dots \tau_k^{n_k}]_a\big) = n_1!\dots n_k!\,\sigma(\tau_1)^{n_1}\dots\sigma(\tau_k)^{n_k},
$$
where $\tau_1,\dots,\tau_k$ are distinct trees with respective multiplicities $n_1,\dots,n_k$. This number does not depend on the label $a$. Denote by $\langle\mathcal{T}\rangle$ the real vector space spanned by $\mathcal{T}$ and by $\langle\mathcal{T}\rangle^*$ its algebraic dual, identified with $\langle\mathcal{T}^*\rangle$, for some copy $\mathcal{T}^*$ of $\mathcal{T}$.

\smallskip

Set $V({\bf 1}^*) := 0$. and define recursively, for $a\in\llbracket 1,\ell\rrbracket$ and $\tau_1,\dots,\tau_n$ in $\mcT$, \textit{vector fields} on $\RR^d$ indexed by dual trees
\begin{equation}
\label{EqDifferentialOperators}
V(\bullet_a^*) := V_a,   \quad   V\big(([\tau_1\dots\tau_n]_a)^*\big) := \frac{1}{\sigma\big({[\tau_1\dots\tau_n]_a}\big)}(D^nV_a)\big(V(\tau_1^*),\dots,V(\tau_n^*)\big);
\end{equation}
\textit{vector fields} on $\RR^d$ are seen as first order differential operators. We define \textit{differential operators} indexed by dual forests setting
\begin{equation}
\label{EqDifferentialOperatorsBis}
V\big(\tau_1^*\cdots\tau_n^*\big)f = \frac{1}{\sigma\big(\tau_1\dots\tau_n\big)}\big(D^nf\big)\big(V(\tau_1^*),\dots,V(\tau_n^*)\big),
\end{equation}
with $\sigma\big(\tau_1\dots\tau_n\big) := \sigma\big({[\tau_1\dots\tau_n]_a}\big)$, for any $1\leq a\leq \ell$. Set $|\tau_1\cdots\tau_n| = |\tau_1| + \cdots + |\tau_n|$, and denote by $\varphi = \tau_1\cdots\tau_n$, a generic forest, and by $\mathcal{F}$ the set of these forests. In those terms, Equation \eqref{EqExpansion} rewrites 
\begin{equation}
\label{EqODEBranched}
f(x_t) = f(x_s) + \sum_{\varphi\in\mathcal{F}; \vert\varphi\vert\leq k_0} H^\varphi_{ts} \, \big(V(\varphi^*) f\big)(x_s) + O\big(|t-s|^{k_0+1}\big).
\end{equation}
The coefficients $H^\varphi_{ts}$ are linear combinations of the initial coefficients $H^I_{ts}$. One has actually the recursive definition
\begin{equation}
\label{EqXTau}
H^{\bullet_a}_{ts} := h^a_t-h^a_s,\qquad H^{[\tau_1\dots\tau_n]_a}_{ts} = \int_s^t \prod_{i=1}^n H^{\tau_i}_{us}\,dh^a_u;
\end{equation}
formulas \eqref{EqODEBranched} and \eqref{EqXTau} go back to Cayley' seminal work \cite{Cayley}. In the same way as the family of iterated integrals of $h$ lives in an algebraic structure, the free nilpotent Lie group, the family of all $H^\tau$ also lives in some algebraic structure, the dual of a Hopf algebra. We refer the reader to Sweedler's little book \cite{Sweedler} for the basics on Hopf algebras; all that we need to know is elementary and recalled below.

\medskip

\textbf{\textsf{2. Hopf algebraic structure --}} Let $(\mcF,\cdot)$ stand for the set of commuting real-valued polynomials with indeterminates the elements of $\mcT$ and polynomial multiplication operation; equip the algebraic tensor product $\mcF\otimes\mcF$ with the induced product $(a\otimes b)(c\otimes d)$ = $(ac)\otimes (bd)$. We define a coproduct $\Delta$ on $\mcF$ as follows. Given a labelled rooted tree $\tau$, denote by $\textsf{Sub}(\tau)$ the set of subtrees of $\tau$  with the same root as $\tau$. Given such a subtree $s$, we obtain a collection $\tau_1,\dots,\tau_n$ of labelled rooted trees by removing $s$ and all the adjacent edges to $s$ from $\tau$. Write $\tau\backslash s$ for the monomial $\tau_1\dots\tau_n$. One defines a linear map $\Delta : \mcT \rightarrow \mcF\otimes\mcF$ by the formula 
$$
\Delta\tau = \sum_{s\in\textsf{Sub}(\tau)}(\tau\backslash s)\otimes s,
$$
and extend it to $\mcF$ by linearity and by requiring that it is multiplicative 
$$
\Delta(\tau_1\dots\tau_n) := \Delta(\tau_1)\dots\Delta(\tau_n). 
$$
This coproduct is coassociative, $(\Delta\otimes\textrm{Id})\Delta = (\textrm{Id}\otimes\Delta)\Delta$, as it should be in any Hopf algebra. Given a real-valued map ${\bf Y}$ on $\mcT$, we extend it into an element of the dual space $\mcF^*$ of $\mcF$ setting 
$$
Y^{\tau_1\dots\tau_n} := \langle {\bf Y}, \tau_1\cdots\tau_n\rangle := \prod_{i=1}^n Y^\tau_i, 
$$
for a monomial $\tau_1\cdots\tau_n$, and by linearity. One defines a convolution product on $\mcF^*$ setting
\begin{equation}
\label{EqConvolution}
({\bf Y\star X})^\tau := {\bf (Y\otimes X)}(\Delta\tau) = \sum_{s\in\textsf{Sub}(\tau)} Y^{\tau\backslash s}X^s,
\end{equation}
for a labelled rooted tree $\tau$. The third ingredient needed to define the Hopf algebra structure of $\mcF$ is an antipode $\frak{S}$, that is a map $\frak{S} : \mcF\rightarrow\mcF$  inverting $\Delta$ in the sense that $\textrm{M}(\textrm{Id}\otimes\frak{S})\Delta = \textrm{M}(\frak{S}\otimes\textrm{Id})\Delta =\textrm{Id}$, where $\textrm{M}$ stands for the multiplication map $\textrm{M}(a\otimes b) = ab\in\mcF$. An explicit formula for $\frak{S}$ was first obtained in \cite{ConnesKreimer}; see \cite{ChartierHairerVilmart} for a simple and enlighting proof. From its definition, the inverse ${\bf a}^{-1}$ of any element $\bf a$ of $\mcF^*$ is given by the formula ${\bf a}^{-1}= \frak{S}^*{\bf a}$, where $\frak{S}^*$ is the dual map in $\mcF^*$ of the antipode map in $\mcF$, that is $\big({\bf a}^{-1},\tau\big) = \big({\bf a},\frak{S}\tau\big)$.

\ssk

Denote  by $\theta=\tau_1\dots\tau_n$, with $\tau_i\in\mcT$, a generic element of the canonical basis $\mcB$ of $\mcF$. Write $({\bf \theta}^*)_{\theta\in\mcB}$ for the dual canonical basis and use bold letters to denote generic elements of $\mcF^*$. For an element $\bf a$ of $\mcF^*$ of the special form ${\bf a} = \sum_{\tau'\in\mcT}a^{\tau'}(\tau')^*$, with $a^{\tau'}\in\bbR$, we have for instance
\begin{equation*}
\begin{split}
{\bf a\star a} &= \sum_{\theta}({\bf a\otimes a})(\Delta\theta)\,\theta^* = \sum_{\tau\in\mcT}({\bf a\otimes a})(\Delta\tau)\,\tau^*  \\
                     &= \sum_{\tau\in\mcT}\Big(\sum a^{\tau_1}a^{\tau_2}\Big)\,\tau^*,
\end{split}
\end{equation*}
where the inside sum is over the set of subtrees $\tau_1$ of $\tau$ with the same root as $\tau$, and such that $\tau\backslash\tau_1$ is a subtree $\tau_2$ of $\tau$. Note in particular that ${\bf a}\star {\bf a}$ is again of the form $\sum_{\tau'\in\mcT} b^{\tau'}(\tau')^*$. More generally, we have for such an ${\bf a}\in\mcF^*$
\begin{equation}
\label{EqAStarN}
{\bf a}^{\star n} = \sum_{\tau\in\mcT}\Big(\sum a^{\tau_1}\cdots a^{\tau_n}\Big)\tau^*,
\end{equation}
where the inside sum is over the set of disjoint rooted subtrees $\tau_1,\dots,\tau_n$ of $\tau$, with respective roots $\alpha_1,\dots,\alpha_n$ such that $\alpha_{i+1}$ is a descendant of $\alpha_i$ in $\tau$, for all $1\leq i\leq n-1$, and any node of $\tau$ is in one of the subtrees $\tau_i$. The convolution product \eqref{EqConvolution} is closely related to the identities \eqref{EqDifferentialOperators} and \eqref{EqDifferentialOperatorsBis} defining the differential operators $V(\tau^*)$ and $V(\tau_1^*\cdots\tau_n^*)$. The following statement is proved by induction on $n+n'$.

\ssk

\begin{lem}
\label{Lemma}
For any dual forests $\tau_1^*\cdots\tau_n^*$ and $\sigma_1^*\cdots\sigma_{n'}^*$, one has
$$
V\big(\tau_1^*\cdots\tau_k^*\big)\,V\big(\sigma_1^*\cdots\sigma_{k'}^*\big) = V\Big((\tau_1^*\cdots\tau_k^*)\star (\sigma_1^*\cdots\sigma_{k'}^*)\Big).
$$
\end{lem}

\medskip

\textbf{\textsf{3. H\"older branched $p$-rough paths --}} We define a Lie bracket on $\mcF^*$ setting 
$$
{\bf [a,b]_\star := a\star b-b\star a},
$$ 
for which the real vector space $\langle\mcT^*\rangle$ spanned by the dual trees is a Lie algebra. The exponential map $\exp_\star : \mcF^* \rightarrow \mcF^*$ and the logarithm map $\log_\star : \mcF^*\rightarrow\mcF^*$ are defined by the usual series
$$
\exp_\star({\bf a}) = \sum_{n\geq 0}\frac{{\bf a}^{\star n}}{n!},\quad \log_\star({\bf b}) = \sum_{n\geq 1}\frac{({\bf 1^*-b})^{\star n}}{n},
$$
with the convention ${\bf a}^{\star 0} = {\bf 1}^*$, where ${\bf 1}^*$ denote the dual of the empty tree. It is clear from the above formula for ${\bf a}^{\star n}$ that we have 
\begin{equation}
\label{EqExpStarTree}
\big(\exp_\star({\bf a}),\tau\big) = \sum_{n=0}^{|\tau|}\frac{1}{n!}\sum a^{\tau_1}\cdots a^{\tau_n},
\end{equation}
with the same inside sum as in \eqref{EqAStarN}. Formula \eqref{EqAStarN} justifies the convergence of the these sums, in the sense that $\big(\exp_\star({\bf a}),\theta\big)$ and $\big(\log_\star({\bf b}),\theta\big)$ are actually finite sums for any monomial $\theta$, provided $({\bf a}, 1) = 0$ -- the pairing $(\cdot,\cdot)$ is a pairing between $\mcF^*$ and $\mcF$. Setting $\mcF^*_0 := \big\{{\bf a}\in\mcF^*\,;\,({\bf a},1)=0 \big\}$ and $\mcF^*_1 := \big\{{\bf a}\in\mcF^*\,;\,({\bf a},1)=1 \big\}$, one can see that $\exp_\star : \mcF^*_0 \rightarrow\mcF^*_1$ and $\log_\star : \mcF^*_1 \rightarrow\mcF^*_0$ are reciprocal bijections; the following result can be seen as a consequence. (See Sweedler's above mentioned book, or Reutenauer's book \cite{Reutenauer}, Theorem 3.2.)

\medskip

\begin{thm}
The pair $\big(\exp_\star(\langle\mcT^*\rangle),\star\big)$ is a group.
\end{thm}

\medskip

Pick an integer $N\geq 0$. Write $\mcF^*_{(k)}$  for the vector space spanned by the monomials $(\tau_1\dots\tau_n)^*$, with $\sum_{i=1}^n|\tau_i| = k$, and denote by $\pi_{\leq N}$ the natural projection from $\mcF^* $ to the quotient space $\mcF^*\backslash \bigoplus_{k\geq N+1}\mcF^*_{(k)}$. Set 
$$
\pi_{\leq N}(\mcF^*) =: \mcF^*_{\leq N}.
$$
The image of $\exp_{\star}(\frak{g})$ by $\pi_N$ is then diffeomorphic to the real vector space $\langle \mcT^*_{\leq N}\rangle$ spanned by $\pi_{\leq }N(\mcT^*)$; it is a Lie group when equipped with the operation $\pi_{\leq N}\circ\star$, which we still denote by $\star$. Write 
$$
\mcG^{N}_\ell\subset \mcF^*_{\leq N}
$$ 
for that Lie group; it plays for branched rough paths the role that the $[p]$-step free nilpotent Lie group plays for weak geometric H{\"o}lder $p$-rough paths. In the same way as Chen's lifting formula \eqref{EqLiftSmooth} gives a geometric to the 

\ssk

Denote by $\bullet_a^*$ the element of the dual canonical basis, dual to $\bullet_a$. It is elementary to check that the $\mcF^*$-valued function ${\frak H}_{ts} = (H^\tau_{ts})_{\tau\in\mcT,\,|\tau|\leq N}$ associated with a smooth $\RR^\ell$-valued path $h$, as defined in \eqref{EqXTau}, considered as a function of the time parameter $t$, satisfies in $\pi_N(\mcF^*)$ the ordinary differential equation
$$
d{\frak H}_{ts} = \sum_{a=1}^\ell {\frak H}_{ts}\star \big(\bullet_a^* dh^a_t\big);
$$
this implies that ${\frak H}_{ts}$ is an element of $\mcG^{(N)}_\ell$ for all times $t\geq s$, and that 
$$
\frak{H}_{ts} = {\bf H}_s^{-1}{\bf H}_t,
$$ with ${\bf H}_r :={\frak H}_{r0}$, for all $r\geq 0$. 

\medskip

\begin{definition*}
Let $2<p$. A \textsf{\textbf{H\"older branched $p$-rough path on $[0,T]$}} is a $\mcG_\ell^{[p]}$-valued path ${\bfX}$ such that 
$$
\underset{0\leq s<t\leq T}{\sup}\,\frac{|X^\varphi_{ts}|}{|t-s|^{\frac{|\varphi|}{p}}}<\infty,
$$ 
for all forests $\varphi=\tau_1\dots\tau_k$ with $|\varphi| := \sum_{i=1}^k|\tau_i|\leq [p]$, and ${\bf X}_{ts} := {\bf X}_s^{-1}\star{\bf X}_t$. 
\end{definition*}

\medskip

The norm of $\bfX$ is defined as
\begin{equation}
\label{DefnNormBranchedRoughPath}
\|{\bfX}\| := \underset{\varphi\in\mcF\,;\,|\varphi|\leq [p]}{\max}\;\underset{0\leq s<t\leq T}{\sup}\,\frac{|X^\varphi_{ts}|}{|t-s|^{\frac{|\varphi|}{p}}};
\end{equation}
we also define a distance $d(\bf X,\bf Y) = \|{\bf X} - {\bf Y}\|$ on the nonlinear set of H\"older branched $p$-rough path. Given $0\leq s\leq t\leq T$ and a H\"older branched $p$-rough path $\bfX$ defined on the time interval $[0,T]$, we define an element of  $\langle \mathcal{T}^*_{\leq [p]} \rangle$ setting 
$$
{\bf \Lambda}_{ts} := \log_\star {\bfX}_{ts};
$$
by definition
$$
\exp_\star({\bf \Lambda}_{ts}) = {\bf X}_{ts} \, \in \mcG^{[p]}_\ell.
$$

\medskip

\textbf{\textsf{4. Rough differential equations driven by H\"older branched $p$-rough paths --}} Let $\bfX$ be a H\"older branched $p$-rough path over $\bR^\ell$, and $\textrm{F}=(V_1,\dots, V_\ell)$ be $\gamma$-H\"older vector fields on $\bbR^d$, for $\gamma>p$. A path $z$ in $\bbR^d$ is said to be a solution path to the rough differential equation
\begin{equation}
\label{EqRDEBranched}
dz_t = \textrm{F}(z_t)\,{d\bfX}_t
\end{equation}
in the sense of Davie if one has 
\begin{equation*}
x_t = x_s + \sum_{\varphi\in\mathcal{F}; \vert\varphi\vert\leq [p]} \big(V({\bfX}_{ts})\textrm{Id}\big)(x_s) + O\big(|t-s|^a\big).
\end{equation*}
for some constant $a>1$. (This is actually Gubinelli's definition \cite{GubinelliBranched} rather than Davie's definition.) The path $z$ is a solution to Equation \eqref{EqRDEBranched} in the sense of Bailleul if one has 
\begin{equation*}
f(x_t) = f(x_s) + \sum_{\varphi\in\mathcal{F}; \vert\varphi\vert\leq [p]} \big(V({\bfX}_{ts}) f\big)(x_s) + O\big(|t-s|^a\big),
\end{equation*}
for some constant $a>1$, for every real-valued function $\gamma$-H\"older function $f$ on $\bbR^d$. (The $O(\cdot)$ term is allowed to depend on $f$.) Theorem \ref{Thm2} states the equivalence of these two notions of definition. Conditions for the existence of a unique solution are given in Gubinelli' seminal paper \cite{GubinelliBranched}.

\medskip

\subsection{Proof of Theorem \ref{Thm2}}
\label{SubsectionProof}

We only need to prove that a solution in Davie' sense is a solution in Bailleul' sense. We proceed as in Section \ref{SectionMain} and prove that the result comes from a Taylor expansion property satisfied by the time $1$ map of an approximate dynamics and the morphism property stated in Lemma \ref{Lemma}. 

\medskip

Recall we write ${\bf \Lambda}_{ts}$ for $\log_\star{\bfX}_{ts}$. Given some times $0\leq s\leq t\leq T$, let $\mu_{ts}$ stand for the well-defined time $1$ map associated with the ordinary differential equation
\begin{equation}
\label{EqApproximateODEBranched}
\dot y_u = V\big({\bf \Lambda}_{ts}\big)(y_u).
\end{equation}
(The assumption that $\gamma>p$ and the vector fields $V_i$ are $\gamma$-Lipschitz ensures that all the vector field $V({\bf \Lambda}_{ts})$ is globally Lipscthiz continuous.)

\ssk

\begin{lem}
\label{LemmaBis}
There is a constant $a>1$ with the following property. For any $\gamma$-H\"older real-valued function $f$ on $\bbR^d$, one has
$$
f\circ \mu_{ts} = f + V\big({\bf X}_{ts}\big)f + O_f\big(|t-s|^a\big),
$$
for a remainder term $O_f\big(|t-s|^a\big)$ that may depend on $f$, and that has a $C^1$-norm bounded above by a constant multiple of $|t-s|^a$.
\end{lem}

\ssk

\begin{Dem}
The Taylor expansion property follows from the morphism property of the $V$-map from Lemma \ref{Lemma}, writing 
\begin{equation*}
\begin{split}
f(y_1) &= f(y_0) + \int_0^1 \big\{V({\bf \Lambda}_{ts})f\big\}(y_{s_1})\,ds_1   \\
 	   &= f(y_0) + \big\{V({\bf \Lambda}_{ts})f\big\}(x) + \int_0^1\int_0^{s_1} \Big\{V({\bf \Lambda}_{ts})V\big({\bf \Lambda}_{ts}\big)f\Big\}(y_{s_2})\,ds_2ds_1   \\
 	   &= f(y_0) + \big\{V({\bf \Lambda}_{ts})f\big\}(x) + \int_0^1\int_0^{s_1} \Big\{V\big({\bf \Lambda}_{ts}^{\star 2}\big)f\Big\}(y_{s_2})\,ds_2ds_1   \\
 	   &= f(y_0) + \big\{V({\bf \Lambda}_{ts})f\big\}(x) + \frac{1}{2}\, \Big\{V\big({\bf \Lambda}_{ts}^{\star 2}\big)f\Big\}(x) + \int_0^1\int_0^{s_1}\int_0^{s_2} \Big\{V({\bf \Lambda}_{ts}^{\star 3})f\Big\}(y_{s_3})\,ds_3ds_2ds_1   \\
\end{split}
\end{equation*}
and, by induction,
\begin{equation*}
\begin{split}
f(y_1) &= f(y_0) + \sum_{k=1}^{[p]} \frac{1}{k!}\,\big\{V\big({\bf \Lambda}_{ts}^{\star k}\big)f\big\}(x) + \int \Big\{V\big({\bf \Lambda}_{ts}^{\star ([p]+1)}\big)f\Big\}(y_{[p]+1}) \,{\bf 1}_{0\leq s_1\leq \cdots\leq s_{[p]+1}} ds   \\
&= \big\{V({\bf X}_{ts})f\big\}(y_0) +  \Theta_f\big(|t-s|^a\big) + \int \Big\{V\big({\bf \Lambda}_{ts}^{\star ([p]+1)}\big)f\Big\}(y_{[p]+1}) \,{\bf 1}_{0\leq s_1\leq \cdots\leq s_{[p]+1}} ds
\end{split}
\end{equation*}
The $\Theta_f\big(|t-s|^a\big)$ term comes from the fact that ${\bf \Lambda}\in \mathcal{H}_N$ and ${\bf X}\in \mathcal{H}_N$ are build from one another using the restriction of the maps $\log_\star$ and $\exp_\star$, so ${\bf X}_{ts}\in\mathcal{H}_N$ coincides with $\sum_{i=0}^N \frac{1}{n!}\,{\bf \Lambda}_{ts}^{\star n}$ up to some terms in $\mathcal{H}\backslash\mathcal{H}_N$ of size $O\big(|t-s|^a\big)$. The conclusion follows then from that explicit representation of the Taylor remainder.
\end{Dem}
 
\ssk
 
It follows from Lemma \ref{LemmaBis} that a path $z$ is a solution in Davie' sense iff 
$$
z_t = \exp\big(V({\bf \Lambda}_{ts})\big)(z_s) + O\big(\vert t-s\vert^{a'}\big),
$$
for some exponent $1<a'\leq \frac{\gamma}{p}$. The end of the proof of Theorem \ref{Thm2} is then identical to the proof of Theorem \ref{Thm}.

\medskip

\noindent \textbf{\textsf{Remarks --}} \textbf{\textsf{1.}} Write $\epsilon^f_{ts}$ for the remainder $f\circ \mu_{ts} - V({\bf X}_{ts})f$. One can see directly that the $\mu_{ts}$ form a $C^1$-approximate flow from the identity
\begin{equation*}
\begin{split}
\mu_{tu}\circ\mu_{us} &= \big\{V({\bf X}_{tu})\textrm{Id}\big\}\circ\mu_{us} + \epsilon^\textrm{Id}_{tu}\circ\mu_{us}   \\
									  &= V({\bf X}_{us}) V({\bf X}_{tu})\textrm{Id} + \epsilon^{V({\bf X}_{tu})\textrm{Id}}
_{us} + \epsilon^{\textrm{Id}_{tu}}\circ\mu_{us}   \\
									  &= V({\bf X}_{ts})\textrm{Id} + \epsilon^{V({\bf X}_{tu})\textrm{Id}}
_{us} + \epsilon^\textrm{Id}_{tu}\circ\mu_{us}   \\
									  &= \mu_{ts} + \epsilon^\textrm{Id}_{ts} + \epsilon^{V({\bf X}_{tu})\textrm{Id}}
_{us} + \epsilon^\textrm{Id}_{tu}\circ\mu_{us}.
\end{split}
\end{equation*}
One then knows from Theorem 1 in \cite{BailleulRMI} that the rough differential equation \eqref{EqRDEBranched} has a unique solution flow, associated with the $C^1$-approximate flow $\mu$. This provides a direct proof that the rough differential equation \eqref{EqRDEBranched} has a unique solution flow.   \vspace{0.15cm}

\textbf{\textsf{2.}} Here as in Section \ref{SectionMain}, working with H\"older rough paths or rough paths controlled by some more general control $w(s,t)$ does not make any difference. Working with branched rough paths over an infinite dimensional space $E$ does not make any difference either, as long as one takes care of working with a symmetric system of cross seminorms generating the topology of $E$, and complete accordingly the different algebraic tensor spaces $E^{\otimes k}$, such as emphasized in \cite{CassWeidner}, Section 3.

\bigskip
\bigskip

\noindent \textcolor{gray}{$\bullet$} {\sf I. Bailleul} - {\small Institut de Recherche Mathematiques de Rennes, 263 Avenue du General Leclerc, 35042 Rennes, France.} ismael.bailleul@univ-rennes1.fr   \vspace{0.3cm}

\end{document}